\documentclass[11pt]{article}
\usepackage{mathrsfs}
\usepackage{amsmath}
\usepackage{amsmath,amssymb,amsfonts,amsthm,fancyhdr}
\usepackage{epsfig,graphicx,picins,picinpar,subfigure}
\usepackage{pstricks}
\usepackage{fancyvrb}
\usepackage[numbers,sort&compress]{natbib}
\begin{document}
\title{\textbf{The concept of primes and the algorithm for counting the greatest common divisor in Ancient China}}
\author{\emph{\textbf{Shaohua Zhang}}}
\date{{\small School of Mathematics, Shandong University,
Jinan,  Shandong, 250100, China\\
E-mail address: shaohuazhang@mail.sdu.edu.cn}}
 \maketitle
\begin{abstract}
When people mention the number theoretical achievements in Ancient
China, the famous Chinese Remainder Theorem always springs to mind.
But, two more of them---the concept of primes and the algorithm for
counting the greatest common divisor, are rarely spoken. Some
scholars even think that Ancient China has not the concept of
primes. The aim of this paper is to show that the concept of primes
in Ancient China can be traced back to the time of Confuciusor
(about 500 B.C.) or more ago. This implies that the concept of
primes in Ancient China is much earlier than the concept of primes
in Euclid's \emph{Elements}(about 300 B.C.) of Ancient Greece. We
also shows that the algorithm for counting the greatest common
divisor in Ancient China is essentially the Euclidean algorithm or
the binary gcd algorithm. Donald E. Knuth said that "the binary gcd
algorithm was discovered by J. Stein in 1961". Nevertheless, Knuth
was wrong. The ancient Chinese algorithm is clearly much earlier
than J. Stein's algorithm.

\vspace{3mm}\noindent \textbf{Keywords:} prime, irreducible numbers,
Chinese Remainder Theorem, the Euclidean algorithm, the binary gcd
algorithm

\vspace{3mm}\noindent \textbf{2000 MR  Subject Classification:}\quad
00A05, 01A05, 11A41, 11A05

\end{abstract}
\section{Introduction}
China was the first nation who studied indeterminate equations. In
the west, indeterminate equations also are called Diophantine
equations. They are indeterminate polynomial equations that allows
the variables to be integers or rational numbers only. The word
Diophantine refers to the Hellenistic mathematician of the 3rd
century, Diophantus of Alexandria, who made a study of such
equations. He was the author of \emph{Arithmetica} and also was one
of the first mathematicians to introduce symbolism into algebra.
However, before the 1st century, the indeterminate equation "wu jia
gong jin" was studied well in Ancient China. The question of "wu jia
gong jin" is just to find the solutions of Diophantine equation
$f=2a+b=3b+c=4c+d=5d+e=6e+a$. Ancient Chinese mathematicians have
found the least positive integral solutions as follows:
$f=721,a=265,b=191,c=148,d=129,e=76$. See \emph{jiu zhang suan shu}
or [1].

\vspace{3mm}When people mention the number theoretical achievements
in Ancient China, the famous Chinese Remainder Theorem always
springs to mind. Undoubtedly, Chinese Remainder Theorem is the
greatest number theoretical-theorem of Ancient China. But, two more
of them---the concept of primes and the algorithm for counting the
greatest common divisor, are rarely spoken. Some scholars even think
that Ancient China has not the concept of primes. The aim of this
paper is to introduce these two great number theoretical
achievements in Ancient China.

\vspace{3mm}We show that  the concept of primes in Ancient China can
be traced back to the time of Confuciusor or more ago. As we know,
Confuciusor was a Chinese thinker and philosopher. His teachings
have deeply influenced the thought and life of many nations. He was
born in 551 B.C. and died in 479 B.C.. Therefore, the time of
Confuciusor is roughly equal to the time of Pythagoras (about 500
B.C.). This implies that the concept of primes in Ancient China is
much earlier than the concept of primes in Euclid's
\emph{Elements}(about 300 B.C.) [2] of Ancient Greece.

\vspace{3mm}We also show that the algorithm for counting the
greatest common divisor in Ancient China is essentially the
Euclidean algorithm [2] or the binary gcd algorithm [3], see Chapter
1 of \emph{Nine Chapters on the Mathematical Art} [1]. Donald E.
Knuth [4] said that "the binary gcd algorithm was discovered by J.
Stein in 1961". Nevertheless, Knuth was wrong. As we know,
"\emph{The Nine Chapters on the Mathematical Art} (pinyin: \emph{jiu
zhang suan shu}) is a Chinese mathematics book, composed by several
generations of scholars from the 10th-2nd century B.C., its latest
stage being from the 1st century A.D.. This book is the one of the
earliest surviving mathematical texts from China, the first being
\emph{suan shu shu} (202 B.C.-186 B.C.) and \emph{zhou bi suan jing
}". Therefore, the ancient Chinese algorithm is clearly much earlier
than J. Stein's algorithm. For the details, see Section 2.

\section{The algorithm for counting the greatest common divisor in Ancient China}
The Euclidean algorithm [2] has several useful variants [5]. Here,
we introduce only the ancient Chinese algorithm for computing the
greatest common divisor of two positive integers, in which all the
divisions by 2 performed can be done using shifts or Boolean
operations. For its details, see Chapter 1 of \emph{Nine Chapters on
the Mathematical Art}[1]. Chapter 1 of \emph{Nine Chapters on the
Mathematical Art}[1] is called "\emph{fang tian}" or \emph{land
surveying}. In this chapter, the reduction of fraction has been
studied. The original text is the following: " ke ban zhe ban zhi,
bu ke ban zhe, fu zhi fen mu, fen zi zhi shu, yi shao jian duo, geng
xiang jian sun, qiu qi deng ye, yi deng shu yue zhi." From this, one
can obtain the following algorithm.

\vspace{3mm}\noindent {\bf  The ancient Chinese algorithm:~~}%
For a positive integer $a$ and a non-negative integer $b$, this
algorithm finds $\gcd(a,b)$. To begin with, we write $a=2^ex$ and
$b=2^fy$ by pre-computing, where $x$ and $y$ are positive odd
integers, $e$ and $f$ are non-negative integers. Note that
$\gcd(a,b)=2^{min\{e,f\}}\gcd(x,y)$. So, finding $\gcd(a,b)$ is
sufficient to find $\gcd(x,y)$.

\textbf{Step 1}. If $x=y$ then output $\gcd(x,y)=x$ and terminate
the algorithm.

\textbf{Step 2}. Compute $x-y=z$, where $z$ is a positive odd
integer. When $y\geq z$, set $x \leftarrow y$ and $y \leftarrow z$.
Otherwise set  $x \leftarrow z$ and  $y \leftarrow y$. And go to
Step 1.

\vspace{3mm}\noindent {\bf A toy example:~~}%
Compute $\gcd(98,63)$.

\vspace{3mm}\noindent {\bf Solution:~~}%
By the ancient Chinese algorithm, we do the following one by one:

98-63=35

63-35=28

35-28=7

28-7=14

14-7=7

7=7

So $\gcd(98,63)=7$.

\vspace{3mm}Clearly, this algorithm is actually Euclid's algorithm
for computing the greatest common divisor of two positive integers
$a$ and $b$ when $a$ and $b$ are not simultaneously even. People
also call the ancient Chinese algorithm \emph{geng xiang jian sun
shu}. Also based on the aforementioned original text, we get the
variant of ancient Chinese algorithm.

\vspace{3mm}\noindent {\bf  The variant of ancient Chinese algorithm:~~}%
For any given positive odd integers $x$ and $y$, without loss of
generality, we assume that $x\geq y$, this algorithm finds their
greatest common divisor $\gcd(x,y)$ as follows.

 \textbf{Step 1}. If $x=y$ then output $\gcd(x,y)=x$ and terminate the algorithm.

 \textbf{Step 2}. Compute $x-y=2^gz$, where $g$ is a non-negative integer, and $z$
is a positive odd integer. When $y\geq z$, set $x \leftarrow y$ and
$y \leftarrow z$. Otherwise set  $x \leftarrow z$ and  $y \leftarrow
y$. And go to Step 1.

\vspace{3mm}\noindent {\bf Another toy example:~~}%
Compute $\gcd(98,63)$.

\vspace{3mm}\noindent {\bf Solution:~~}%
By the aforementioned variant of ancient Chinese algorithm, we do
the following one by one:

98-63=35

63-35=28

28=4$\times$7

35-7=28

28=4$\times$7

7=7

So $\gcd(98,63)=7$.

\vspace{3mm} By the aforementioned algorithms, one can write the
binary gcd algorithm as follows [6].

\vspace{3mm}\noindent {\bf  The binary gcd algorithm:~~}%
Given two non-negative integers $a$ and $b$, this algorithm finds
their gcd.

\textbf{Step 1}. If $a<b$ exchange  $a$ and $b$. Now, if $b=0$,
output $\gcd(a,b)=a$ and terminate the algorithm. Otherwise, set
$r\leftarrow a\mod b$, $a\leftarrow b$ and $b\leftarrow r$.

\textbf{Step 2}. If $b=0$ output $\gcd(a,b)=a$ and terminate the
algorithm. Otherwise, set $k\leftarrow 0$, and then while $a$ and
$b$ are both even, set $k\leftarrow k+1$, $a\leftarrow a/2$,
$b\leftarrow b/2$.

\textbf{Step 3}. If $a$ is even, repeat $a\leftarrow a/2$ until $a$
is odd. Otherwise, if $b$ is even, repeat $b\leftarrow b/2$ until
$b$ is odd.

\textbf{Step 4}. Set $t\leftarrow (a-b)/2$. If $t=0$, output
$\gcd(a,b)=2^ka$ and terminate the algorithm.

\textbf{Step 5}. While $t$ is even, set $t\leftarrow t/2$. Then if
$t>0$ set $a\leftarrow t$, else set $b\leftarrow -t$ and go to Step
4.

\vspace{3mm}\noindent {\bf Remark 1:~~}%
A challenging mathematical problem is to find an asymptotic estimate
for the number of steps and the number of shifts performed in the
binary gcd algorithm. See [4] and [6].

\section{The concept of primes in Ancient China}
From Euclid's \emph{Elements} [2], we know that the concept of
primes is one of the important number theoretical achievements in
Ancient Greece. Generally speaking, the concept of primes in Ancient
China goes back to Shanlan Li. In [7], Dunjie Yan  pointed out that
Shanlan Li proved that Fermat's little theorem and gave four methods
for testing whether a number is prime or not. This leads that some
scholars even think that Ancient China has not the concept of primes
because Shanlan Li was born in 1811 and died in 1882. Maybe, the
main reason that people have been misguided is of that the work of
Ancient Chinese mathematicians has not been disseminated well.

\vspace{3mm} Shanlan Li was not the first person in China who has
the concept of primes at all. By internet, we learn that Guoping
Kong pointed out that Hui Yang was the first to propose the concept
of primes in China. In 1274, Hui Yang published his book \emph{Cheng
Chu Tong Bian Ben Mo } which means Alpha and omega of variations on
multiplication and division appeared. In this book, Hui Yang gave a
quick algorithm for multiplication. He said:"cheng wei fan zhe, yue
wei er duan, zuo er ci cheng zhi, shu ji wei jian er yi cheng, zi ke
wu wu ye." What he said is of that in order to count $a\times b$,
one had better factor firstly $a$ or $b$ such that $a=c\times d$ or
$b=e\times f$, additionally, count $a\times e=m$, and $m\times f=n$,
and so on, finally, $a\times b=n$. For example, let's compute
$38367\times23121$. Firstly, we factor $23121=3^2\times7\times 367$.
Secondly, compute $38367\times9=345303$, $345303\times7=2417121$ and
$2417121\times 367=887083407$. So, $38367\times23121=887083407$.
Clearly, Hui Yang's method looks efficient once one can factor
quickly the multiplicand or multiplier. But, as we know, factoring a
large number $n$ indeed is hard unless $n$ has only small prime
divisors. Fortunately, due to the requirement of his method, Hui
Yang gave the definition of irreducible numbers. A number $n>1$ is
irreducible if $n$ can not be factored as the product of two numbers
which are less than $n$. Hui Yang's irreducible numbers in fact are
prime numbers. It enables us to understand that there are different
understandings of "primes" in Ancient China. Furthermore, he listed
all irreducible numbers (primes) between 200 and 300 as follows:
211, 223, 227, 229, 233, 239, 241, 251, 257, 263, 269, 271, 277,
281, 283 and 293. Hui Yang's work on irreducible numbers is very
significative.

\vspace{3mm} Was Hui Yang the first person in Ancient China who has
the concept of primes? Surely, the answer is no. Clearly, the
concept of "irreducible" is very closely related the reduction of
fraction. And the reduction of fraction has been studied well in
\emph{Nine Chapters on the Mathematical Art}.  It's also worth
noting that Hui Yang studied carefully \emph{Nine Chapters on the
Mathematical Art}, and in 1261, he wrote the \emph{Xiang jie jiu
zhang suan fa} (Detailed analysis of the mathematical rules in the
\emph{Nine Chapters on the Mathematical Art} and their
reclassifications). If the concept of "irreducible" is very closely
related to the reduction of fraction, then it should go back to
\emph{suan shu shu} (202 B.C.-186 B.C.) because in this book, the
reduction of fraction has been studied well. From this, one can see
again that Ancient Chinese mathematicians (202 B.C.-186 B.C.) had
algorithms for counting the greatest common divisor of two integers.
Did \emph{Nine Chapters on the Mathematical Art} grow out of
\emph{suan shu shu}? I think that is reasonable, but not sure.
Anyway, the concept of primes in Ancient China should exist prior to
\emph{suan shu shu} which perhaps is the first mathematical book of
China.

\vspace{3mm}Moreover, from the preface of his magnum opus
\emph{History of the Theory of Numbers, Volume I: Divisibility and
Primality}, Leonard Eugene Dickson [8] wrote: "Fermat stated in 1640
that he had a proof of the fact, now known as Fermat's little
theorem, that, if $p$ is any prime and $x$ is any integer not
divisible by $p$, then $x^{p-1}-1$ is divisible by $p$. This is one
of the fundamental theorem of the theory of numbers. The case $x=2$
was known to the Chinese as early as 500 B.C."  From this, it can be
seen that there must be the concept of primes in Ancient China. In
the first paragraph of the third chapter of [8], Dickson wrote
again: "The Chinese seem to have known as early as 500 B.C. that
$2^p-2$ is divisible by the prime $p$." Then, he wrote [8, pp91]:
"In a Chinese manuscript dating from the time of Confucius it is
stated erroneously that $2^{n-1}-1$ is not divisible by $n$ if $n$
is not prime." This is called the "Chinese hypothesis" or "Chinese
congruence" today. As we mentioned in Section 1, the time of
Confuciusor is roughly equal to the time of Pythagoras (about 500
B.C.). Although Ancient Chinese mathematicians perhaps have not
found that $2^{341-1}-1$ is divisible by $341=11\times31$, we see
again that the ancient Chinese as early as 500 B.C must have the
concept of primes from this. This implies that the concept of primes
in Ancient China is much earlier than the concept of primes in
Euclid's \emph{Elements}(about 300 B.C.) of Ancient Greece.

\vspace{3mm}\noindent {\bf Remark 2:~~}%
In his books \emph{The New Book of Prime Number Records} or
\emph{The little book of the bigger primes}, Ribenboim, P. pointed
out "it is incorrect to ascribe this question (the Chinese
congruence) to the ancient Chinese......" However, in their book
\emph{su shu pan ding yu da shu fen jie}(\emph{The determination of
primes and decomposition of large numbers}), Qi Sun and Jinghua
Kuang pointed out that the Chinese congruence has been studied in
\emph{Book of Changes} of Ancient China. As we know, \emph{Book of
Changes} is also called \emph{Zhou yi}, \emph{ Classic of Changes},
\emph{I Ching}, \emph{Yi King} or \emph{Yi Jing} (Pinyin). It is one
of the oldest of the Chinese classic texts more than 3000 years ago.
It not only is  a philosophical and divinatory book but also has
mathematical significance. It is the source of the binary numeral
system. "Richard S. Cook reported that \emph{Book of Changes}
demonstrated a relation between the golden ratio and linear
recurrence sequences." And so on. G.W. Leibniz studied carefully
\emph{Book of Changes}and believed incorrectly that the Chinese
congruence holds. The aim of this paper is not to talk about the
Chinese congruence. From this, we believe again that there must be
the concept of primes in Ancient China although there are different
understandings of "primes".

\clearpage
\end{document}